\documentclass{elsart}

\usepackage{amssymb,amsmath,graphicx}
\usepackage{tikz}

\date{}
\def\NN{\hbox{\sf I\kern-.13em\hbox{N}}}
\def\RR{\hbox{\sf I\kern-.14em\hbox{R}}}
\def\ZZ{{\hbox{\sf Z\kern-.43emZ}}}

\def\Aut {\mathop{\rm Aut}\nolimits}

\def\mod {\mathop{\rm mod}\nolimits}

\newtheorem{theorem}{Theorem}[section]

\newtheorem{lemma}[theorem]{Lemma}

\newtheorem{corollary}[theorem]{Corollary}

\begin{document}
 
 \begin{frontmatter}
 
\title{Optimal $L(d,1)$-labeling of certain direct graph bundles cycles over cycles and Cartesian graph bundles cycles over cycles }

\author{Irena Hrastnik Ladinek} 

\address{
Faculty of Mechanical Engineering, University of Maribor, Smetanova 17,\\ 2000 Maribor, Slovenia; irena.hrastnik@um.si } 

\begin{abstract}
An $L(d,1)$-labeling of a graph $G$ is an assignment of nonnegative integers to the vertices such that adjacent vertices receive labels that differ by at least $d$ and those at a distance of two receive labels that differ by at least one, where $d\geq 1$. Let $\lambda^d_1 (G)$ denote the least $\lambda$ such that $G$ admits an $L(d,1)$-labeling using labels from $\{0,1,\ldots , \lambda \}$. We prove that $\lambda^d_1(X)\leq 2d+2$ for certain direct graph bundle $X= C_m\times^{\sigma_\ell} C_n$ and certain Cartesian graph bundle $X= C_m\Box^{\sigma_\ell} C_n$, where $\sigma_\ell$ is a cyclic $\ell$-shift, with equality if $1\leq d\leq 4$.
\end{abstract}

\begin{keyword} $L(d,1)$-labeling, $\lambda$-number, direct product of graph, direct graph bundle, Cartesian product of graph, Cartesian graph bundle, cyclic $\ell$-shift, channel assignment.
\end{keyword}


 \end{frontmatter} 

\section{Introduction}

The Frequency Assignment Problem (FAP) is a combinatorial optimization problem that arises in the field of telecommunications and radio frequency (RF) engineering. The goal of the FAP is to assign a set of communication frequencies to a set of transmitters while satisfying certain constraints and minimizing interference. This problem was first formulated as a graph coloring problem in 1980 by Hale \cite{Hale}. By Roberts \cite{Roberts}, in order to avoid interference, any two ''close'' transitters must receive different channels and any two ''very close'' transmitters must receive channels that are at least two channels apart. To translate the problem into the language of graph theory, the transmitters are represented by the vertices of a graph; two vertices are ''very close'' if they are adjacent and ''close'' if they are of distance two in the graph. Based on this problem,  Griggs and Yeh \cite{GrYe} introduced $L(2,1)$-labeling on a simple graph. 

Formally, an $L(d,1)$-labeling of a graph $G$ is an assignment $f$ of non-negative integers to vertices of $G$ such that $$|f(u)-f(v)|\geq \left\{\begin{array}{lr}d;& d(u,v)=1,\\
1;& d(u,v)=2,
\end{array}
\right.$$  where $d\geq 1$. The difference between the largest label and the smallest label assigned by $f$ is called the span of $f$, and the minimum span over all $L(d,1)$-labeling of $G$ is called the $\lambda_1^d$-number of $G$, denoted by $\lambda_1^d(G)$. The general problem of determining $\lambda_1^d(G)$ is $NP$-hard \cite{GeoMau}.

The $L(2,1)$-labeling has been extensively studied in recent past by many researchers, the common trend is either to determine the value of $L(2,1)$-labeling number or to suggest bounds for particular classes of graphs \cite{ChaKuo,Gonc,KlaVes,Kral,Whit}. 

Graph products are one of the natural constructions giving more complex graphs from simple ones.
Graph bundles \cite{Pisanski,PisanskiVra}, also called twisted products, are a generalization of product graphs, 
which have been  (under various names) frequently used as
computer topologies or communication networks.
A famous example is the ILIAC IV supercomputer \cite{ILLIAC}.
While the labeling problem of Cartesian and direct products are well studied \cite{Chiang,Jha,JhaKla,Klav,Schwarz,Shao}, 
there is much less known on the labeling problem of Cartesian  and direct graph bundles \cite{KlavMo,Pisanski}.

We shall need the following lemma, see \cite{GeMa}.

\begin{lemma}
\label{lemma 1}
If $G$ is a graph with maximum degree $\Delta$  and $G$ includes a vertex with $\Delta$ neighbours, each of which is of degree $\Delta$, then  $\lambda^d_1(G)\geq \Delta+2d-2$, where $1\leq d\leq \Delta$. \;\; \rule{2mm}{2mm}
\end{lemma} 

The central result of this paper is that for certain direct graph bundle $C_m\times^{\sigma_\ell} C_n$ and certain Cartesian graph bundle $C_m\Box^{\sigma_\ell} C_n$, where $\sigma_\ell$ is a cyclic $\ell$-shift, the preceding lower bound corresponds to the exact value of $\lambda_1^d$. Analogous result is known with respect to $\lambda_1^d$-numbering of direct products of cycles and Cartesian product of cycles \cite{Jha}.

The rest of the paper is organized as follows. 
In the next section, we provide the basic definitions and some preliminary observations that are needed for the outline of our results. Section 3 deals with the $\lambda_1^d$-numbering of direct graph bundles cycles over cycles while Section 4 presents the analogous result for Cartesian graph bundles cycles over cycles. Methods of attack are similar.

\section{Terminology and Preliminaries} 

A finite, simple and undirected graph $G=(V(G),E(G))$ is given by a set of vertices $V(G)$ and a set of edges $E(G)$. As usual, the edge $\{i,j\}\in E(G)$ is shortly denoted by $ij$. Although here we are interested in undirected graphs, the order of the vertices will sometimes be important, 
 for example when we will assign  automorphisms to edges of the base graph.
 Is such case we assign two opposite arcs $\{(i,j), (j,i)\}$ to edge $\{i,j\}$. 
  
Two graphs $G$ and $H$ are called {\em isomorphic}, in symbols $G\simeq H$, if there exists a bijection $\varphi $ 
from $V(G)$ onto $V(H)$ that preserves adjacency and nonadjacency.
In other words, a bijection $\varphi : V(G) \to v(H)$  is an  {\em isomorphism} when:   $\varphi(i)\varphi(j)\in E(H)$ 
if and only if $ij \in E(G)$.  An isomorphism of a graph $G$ onto itself is called an {\em automorphism}.
The identity automorphism on $G$ will be denoted by $id_G$ or shortly $id$.

The {\em cycle} $C_n$ on $n$ vertices is defined by  $V(C_n) = \{0,1,\dots,n-1\}$ and $ij \in E(C_n)$ if $i=(j\pm 1) \mod n$.
Denote by $P_n$ the {\em path} on $n\geq 1$ distinct vertices $0,1,2,\ldots, n-1$ with edges $ij \in E(P_n)$ if $j=i+1, 0\leq i <n-1.$
For a graph $G=(V,E)$ the {\em distance} $d_G(u,v)$, or briefly $d(u,v)$, between vertices $u$ and $v$ is defined as the number of edges on a shortest $u,v$-path.

Automorphisms of a cycle are of two types: {\em cyclic shift of the cycle by $\ell$ elements} will be briefly called cyclic $\ell$-shift and {\em reflections with one, two or no fixed points} (depending on parity of $n$). We will focus on the first type. A cyclic $\ell$-shift, denoted by $\sigma_\ell$, $0\leq \ell <n$, defined as $\sigma_\ell(i)=(i+\ell)\mod n$ for $i=0,1,\ldots,n-1$. As a special case we have the identity ($\ell=0$).

Let $G=(V(G),E(G))$ and $H=(V(H),E(H))$ be connected graphs. The {\em  direct product} $G\times H$ and the {\em Cartesian
product} $G\Box H$ of $G$ and $H$ are defined as follows: $V (G \times H ) = V (G\Box H ) = V(G) \times V(H)$;
$E(G\times H)=\{\{(g_1,h_1)(g_2,h_2)\}~|~ \{g_1,g_2\}\in E(G) {\mathop{\rm \;and\;}} \{h_1,h_2\}\in E(H)\}$ and $E(G\Box H)=\{\{(g_1,h_1)(g_2,h_2)\}~|~ \{g_1,g_2\}\in E(G) {\mathop{\rm \;and\;}} h_1=h_2, {\mathop{\rm \;or\;}} \{h_1,h_2\}\in E(H) {\mathop{\rm \;and\;}} g_1=g_2\}$.
For more facts on the direct and the Cartesian product of graphs we refer to \cite{ImKl}. 

Let $B$ and $G$ be graphs and $\Aut(G)$ be the set of automorphisms of $G$.
To any ordered pair of adjacent vertices $u,v \in V(B)$  we will assign an automorphism of $G$. 
Formally, let $\sigma : V(B)\times V(B) \to \Aut(G)$.
For brevity, we will write $\sigma(u,v)  = \sigma_{u,v}$
and assume  that $\sigma_{v,u} = \sigma_{u,v}^{-1}$ for any $u,v$ $\in V(B)$.

Now we construct the graphs $X_1$ and $X_2$ as follows.
The vertex set of $X_1$ and $X_2$ is the Cartesian product of vertex sets, $V(X_1) =  V(X_2) =V(B) \times V(G)$.
The edges of $X_1$ are given by the rule:
for any $b_1b_2 \in E(B)$ and any $g_1g_2 \in E(G)$, the vertices $(b_1,g_1)$ and 
$(b_2,\sigma_{b_1,b_2}(g_2)) $ are adjacent  in $X$.
We call $X_1$ a {\em direct graph bundle}  with base   $B$ and fibre $G$ and write $X_1 = B\times^\sigma G$.
The edges of $X_2$ are given by the rule:
for any  $g_1g_2 \in E(G)$ and any $b \in V(B)$, the vertices $(b,g_1)$ and $(b,g_2) $ are adjacent in $X_2$, and for any $b_1b_2 \in E(B)$ and any  $g \in V(G)$, the vertices $(b_1,g)$ and $(b_2,\sigma_{b_1,b_2}(g)) $ are adjacent in $X_2$. 
We call $X_2$ a {\em Cartesian graph bundle}  with base   $B$ and fibre $G$ and write $X_2 = B\Box^\sigma G$.

Clearly, if all $\sigma_{u,v}$ are identity automorphisms, the direct graph bundle is isomorphic to 
the direct product  $B\times^\sigma G = B\times  G$  and the Cartesian graph bundle to the Cartesian product  $B\Box^{\sigma} G = B\Box  G$.

A graph bundle over a cycle can always be constructed in a way that all but at most one automorphism are identities.
Fixing $V(C_n)= \{0,1,2,\dots,n-1\}$, let us denote $\sigma_{n-1,0} = \alpha$, 
$\sigma_{i-1,i} = id$ for $i=1,2,\dots,n-1$ and write   
$ C_n\times^\alpha G \simeq C_n\times^\sigma G $ and  $ C_n\Box^\alpha G \simeq C_n\Box^\sigma G $. In this article we will frequently use this fact.

A graph bundle $ C_n\times^\alpha G$ can also be represented as the graph obtained from the product $P_n\times G$ by adding  a copy of $K_2\times G$ between vertex sets $\{n-1\}\times V(G)$ and $\{0\}\times V(G)$ such that if $V(K_2)=\{1,2\}$ and $(1,u)$ is adjacent to $(2,v)$ in $K_2\times G$, then $(n-1,u)$ and $(0,\alpha(v))$ are connected by an edge in $ C_n\times^\alpha G$. In a similar way, graph bundle $C_n\Box^\alpha G$ can be represented from the product $P_n\Box G$ by adding  a copy of $K_2\Box G$ between vertex sets $\{n-1\}\times V(G)$ and $\{0\}\times V(G)$.

The following facts will be useful in the sequel.

\begin{claim}
\label{claim 1}
If $a,b$ and $n$ are integers with $n\geq 1$, then $|(a\mod n)-(b \mod n)|=(|a-b|\mod n)$ or $n-(|a-b|\mod n)$. \;\; \rule{2mm}{2mm}
\end{claim}

\begin{corollary}
\label{posledica 1}
If $a,b,d$ and $n$ are integers with $d,n\geq 1$, then $|(a\mod n)-(b \mod n)|\geq d\Leftrightarrow d\leq  (|a-b|\mod n) \leq n-d$. \;\; \rule{2mm}{2mm}
\end{corollary}

\begin{corollary}
\label{posledica 2}
If $a,b$ and $n$ are integers with $n\geq 1$, then $|an-b| \mod n=(b \mod n)$ or $n-(b \mod n)$. \;\; \rule{2mm}{2mm}
\end{corollary}


\section{Direct graph bundle $C_m\times^{\sigma_\ell} C_n$  }

\begin{theorem}
\label{izrek1}
Let $d\geq 1$, $m\geq 3$ and $s=2d+3$. Let $X= C_m\times^{\sigma_\ell} C_n$ be a direct graph bundle with fibre $C_n$ and base $C_m$.  If $n$ is a multiple of $s$ and $\ell$ has a form of $\ell =[ks+(-1)^a 2m]\mod n$ or of $\ell =[ks-(-1)^a (d+1)m]\mod n$, where $a\in \{1,2\}$ and $k\in \ZZ$, then $\lambda^d_1(X)\leq 2d+2$, with equality if $1\leq d\leq 4$.
\end{theorem}

\noindent {\bf Proof:}

To prove this theorem, we will present several  labelings of $X$. We will define four  $L(d,1)$-labelings of $X$ using labels $0,1,\ldots , 2d+2$ according to the cyclic $\ell$-shift $\sigma_\ell$. Let $v=(i,j)\in V(X)$.
\begin{enumerate}  
\item Let $\ell =[ks+(-1)^a2m]\mod n$, where $a\in \{1,2\}$ and $k\in \ZZ$. Define labeling $f_a$ of $v$ as $$f_a(v)=\left[i+(d+a)j\right]\mod s$$ 
\item Let $\ell =[ks-(-1)^a(d+1)m]\mod n$, where $a\in \{1,2\}$ and $k\in \ZZ$. Define labeling $g_a$ of $v$ as $$g_a(v)=\left[(d+a)i+j\right]\mod s$$
\end{enumerate}

All assignments are clearly well-defined. 

Since $ C_m\times^{\sigma_\ell} C_n$ can be represented  as the graph obtained from the product $P_m\times C_n$ by adding edges between vertex sets $\{m-1\}\times V(C_n)$ and $\{0\}\times V(C_n)$ (i.e. edges corresponding to the nontrivial automorphism ${\sigma_\ell} ,\ell\ne 0$) first observe the product $P_m\times C_n$. 

First, let $X$ have labeling $f_a$. Then $v$ be assigned the integer  $$f_a(v)=\left[i+(d+a)j\right]\mod s.$$ 
Let $w$ be a vertex adjacent to $v$ in $P_m\times C_n$, so $w$ is of the form $(i+i', (j+j')\mod n)$ with one of the following properties:
\begin{enumerate}
\item  $i=0, i'=1,|j'|=1$,
\item  $i=m-1, i'=-1,|j'|=1$, 
\item $i=1,2,\ldots m-2, |i'|=1,|j'|=1$.
\end{enumerate}

It is clear that $$f_a(w)=\left[i+(d+a)j+i'+(d+a)j'\right]\mod s.$$
By Corollary \ref{posledica 1}, to show that $|f_a(v)-f_a(w)|\geq d$, it is enough to show that $$d\leq \left[|i'+(d+a)j'|\mod s\right] \leq s-d=d+3$$  and, by symmetry, if $X$ has labeling $g_a$  that $$d\leq \left[|(d+a)i'+j'|\mod s\right] \leq d+3.$$ The reader can easily verify that $|i'+(d+a)j'|\mod s,|(d+a)i'+j'|\mod s\in \{d,d+1,d+2,d+3\}$.

Let now $z$ be a vertex at a distance of two from $v$ in $P_m\times C_n$. Then $z$ is of the form $(i+i'', (j+j'')\mod n)$ with one of the following properties:

\begin{enumerate}
\item $i\in \{0,1\}, i''\in \{0,2\},j''\in \{2,0,-2\}$, $i''$ and $j''$ are not all zero,
\item $i\in \{m-2,m-1\}, i''\in \{0,-2\},j''\in \{2,0,-2\}$, $i''$ and $j''$ are not all zero,
\item $i\in \{2,3,\ldots ,m-3\}, i'',j''\in \{2,0,-2\}$, $i''$ and $j''$ are not all zero.
 
\end{enumerate}

Let $X$ have labeling $f_a$. Note that $z$ receives the label $$f_a(z)=\left[i+(d+a)j+i''+(d+a)j''\right]\mod s.$$ 
To show that $|f_a(v)-f_a(z)|\geq 1$, it is enough to show, by Collorarly \ref{posledica 1}, that $$1\leq \left[|i''+(d+a)j''|\mod s\right] \leq 2d+2 $$ and, 
if $X$ has labeling $g_a$, that $$1\leq \left[|(d+a)i''+j''|\mod s\right] \leq 2d+2 .$$
Since  $|i''+(d+a)j''|\mod s,|(d+a)i''+j''|\mod s\in \{1,2, 3, 2d, 2d+1, 2d+2\}$, it follows. 


%

We will now consider the edges between the fiber over $m-1$ and the fiber over $0$ in $X= C_m\times^{\sigma_\ell} C_n$. These are edges $(m-1,u)(0,(u+j'+\ell)\mod n)$, where $u\in V(C_n)$ and $|j'|=1$ or in another form $(0,u)(m-1,(u+j'-\ell)\mod n)$, where $ u\in V(C_n)$ and $|j'|=1$ (recall that $\sigma_{\ell}(i)=(i+\ell)\mod n$ and $\sigma_{\ell}^{-1}(i)=(i-\ell)\mod n$).


First, let's consider two adjacent vertices, one from the fiber over $m-1$ and the other from the fiber over $0$.
Let $X$ has labeling $f_a$ and let $v=(m-1,j)$. If $w$ is a neighbour from the fiber over $0$, then $w$ is of the form $w=(0,(j+j'+\ell)\mod n)$, where $|j'|=1$ and  $\ell =[ks+(-1)^a2m]\mod n$ for some $k\in \ZZ$. Then  $$f_a(v)=\left[m-1+(d+a)j\right]\mod s$$ and $$f_a(w)=\left[(d+a)(j+j'+\ell)\right]\mod s.$$ To show that $|f_a(v)-f_a(w)|\geq d$, it is enough to show that 
\begin{equation}
\label{a}
d\leq [|m-1-(d+a)(j'+\ell)|\mod s]\leq d+3,
\end{equation} by Corollary  \ref{posledica 1}.

Because $\ell =[ks+(-1)^a2m]\mod n$ for some $k\in \ZZ$, there exists $k'\in \ZZ$ such that $\ell =k's+(-1)^a2m$. Hence $$|m-1-(d+a)(j'+\ell)|\mod s=|m-1-(d+a)(j'+k's+(-1)^a2m)|\mod s$$ 
$$=|m(1-(-1)^a2d-(-1)^a2a)-1-(d+a)j'-(d+a)k's)|\mod s $$
For $a=1$ we get $$|ms-1-(d+1)j'-(d+1)k's)|\mod s =$$ $$|s(m-(d+1)k')-(1+(d+1)j')|\mod s$$ and 
for $a=2$  $$|-ms-1-(d+2)j'-(d+2)k's)|\mod s =$$ $$|s(-m-(d+2)k')-(1+(d+2)j')|\mod s.$$
Since $(1+(d+a)j')\mod s\in \{d+2,d+3\}$ and $s-(1+(d+a)j')\mod s\in  \{d,d+1\},$ the claim is true, by Corollary \ref{posledica 2}.

Let now $X$ have labeling $g_a$. Then $\ell =k's-(-1)^a(d+1)m$ for some $k'\in \ZZ$ and for adjacent vertices $v=(m-1,j)$ and $w=(0,(j+j'+\ell)\mod n)$ we need to show that 
\begin{equation}
\label{b}
d\leq [|(d+a)(m-1)-(j'+\ell)|\mod s]\leq d+3.
\end{equation}

Note that $$|(d+a)(m-1)-(j'+\ell)|\mod s=|(d+a)(m-1)-(j'+k's-(-1)^a(d+1)m)|\mod s$$
$$=|m(d+a+(-1)^ad+(-1)^a)-(d+a)-j'-k's|\mod s.$$
For $a=1$ we get
$$|-k's-((d+1)+j')|\mod s$$ and for $a=2$ $$|ms-(d+2)-j'-k's|\mod s=|s(m-k')-((d+2)+j')|\mod s.$$

Since $((d+1)+j')\mod s,s-((d+2)+j')\mod s\in\{d,d+2\}$ and $((d+2)+j')\mod s,s-((d+1)+j')\mod s\in\{d+1,d+3\}$, the desired result follows.

Finally, observe vertices at a distance of two where the shortest path between them contains at least one edge between the fiber over $m-1$ and the fiber over $0$. Let $v$ and $z$ be  two such vertices. We claim that the label of $v$ is not equal to the label of $z$.

First, let $v$ and $z$ be vertices from the  fibre over $m-1$ (it is analogous if both vertices are from the fibre over 0). Let $v=(m-1,j)$. Then $z$ is of the form  $z=(m-1,(j+j')\mod n)$, where $|j'|=2$ (using the fact that $(m-1,j)$ and $(m-1,(j+2) \mod  n)$ have a common neighbour $(0,(j+1+\ell)\mod n)$
and  $(m-1,j)$ and $(m-1,(j-2) \mod  n)$ have a common neighbour $(0,(j-1+\ell)\mod n)$). We have already considered such vertices when we considered vertices at a distance  of two in the graph $P_m\times C_n$.  

Let now $v$ be a vertex from the fiber over $m-1$ and $z$ a vertex from the fiber over  $1$ (similarly, $v$ is from the fibre over $m-2$ and $z$ is from the fiber over $0$). 

Let $v=(m-1,j)$. Then $z$ is of the form $z=(1,(j+j'+\ell)\mod n)$, where $j'\in \{-2,0,2\}$.  Let $X$ have labeling $f_a$: $$f_a(v)=\left[m-1+(d+a)j)\right]\mod s$$ and $$f_a(z)=\left[1+(d+a)(j+j'+\ell)\right]\mod s.$$ 

In this case is enough to show, by Corollary \ref{posledica 1}, that  $$1\leq [|m-2-(d+a)(j'+\ell)|\mod s]\leq 2d+2.$$ 

In the proof of (\ref{a}), we showed that $m-1-(d+a)(j'+\ell)=k''s-(1+(d+a)j')$ for some $k''\in \ZZ$. Therefore $$|m-2-(d+a)(j'+\ell)|\mod s=|k''s-(2+(d+a)j')|\mod s.$$ Since $(2+(d+a)j')\mod s\in \{1,2,3\}$ and $s-(2+(d+a)j')\mod s\in \{2d,2d+1,2d+2\},$  the claim is true, by Corollary \ref{posledica 2}.\\

Now  observe labeling $g_a$. We will consider similar to above. We claim that $$1\leq [|(d+a)(m-2)-(j'+\ell)|\mod s]\leq 2d+2.$$ In the proof of (\ref{b}) we see that $(d+a)(m-1)-(j'+\ell)=k''s-((d+a)+j')$ for some $k''\in \ZZ$. Therefore $$|(d+a)(m-2)-(j'+\ell)|\mod s=|k''s-(2(d+a)+j')|\mod s.$$ Since $(2(d+a)+j')\mod s ,s-((2(d+a)+j')\mod s)\in \{1,3,2d,2d+2\}$, the desired result follows. Accordingly, two vertices that are at a distance of two from each other receive different  labels.\\

We showed that $\lambda^d_1(C_m\times^{\sigma_\ell} C_n)\leq 2d+2$. Further $C_m\times^{\sigma_\ell} C_n$ being a regular graph of degree 4, an application of Lemma \ref{lemma 1} to the preceding statement shows that $\lambda^d_1(C_m\times^{\sigma_\ell} C_n)= 2d+2$, if $1\leq d \leq 4$.

\rule{2mm}{2mm}

The foregoing scheme is illustrated in Fig. \ref{slika1}, where an $L(2,1)$-labeling of $P_9\times P_7$ appears toward that of  $C_9\times^{\sigma_\ell}C_7$ for $\ell=1,3,4,6$.



\begin{figure}[h]
\label{slika1}
\begin{center}
\begin{tikzpicture} [scale=0.65]

	{    				
\filldraw (1,1) circle (2pt); \filldraw (2,1) circle (2pt); \filldraw (3,1) circle (2pt); \filldraw (4,1) circle (2pt);
\filldraw (5,1) circle (2pt); \filldraw (6,1) circle (2pt); \filldraw (7,1) circle (2pt); \filldraw (8,1) circle (2pt); \filldraw (9,1) circle (2pt);
\filldraw (1,2) circle (2pt); \filldraw (2,2) circle (2pt); \filldraw (3,2) circle (2pt); \filldraw (4,2) circle (2pt);
\filldraw (5,2) circle (2pt); \filldraw (6,2) circle (2pt); \filldraw (7,2) circle (2pt); \filldraw (8,2) circle (2pt); \filldraw (9,2) circle (2pt);
\filldraw (1,3) circle (2pt); \filldraw (2,3) circle (2pt); \filldraw (3,3) circle (2pt); \filldraw (4,3) circle (2pt); 
\filldraw (5,3) circle (2pt); \filldraw (6,3) circle (2pt); \filldraw (7,3) circle (2pt); \filldraw (8,3) circle (2pt); \filldraw (9,3) circle (2pt);
\filldraw (1,4) circle (2pt); \filldraw (2,4) circle (2pt); \filldraw (3,4) circle (2pt); \filldraw (4,4) circle (2pt);
\filldraw (5,4) circle (2pt); \filldraw (6,4) circle (2pt); \filldraw (7,4) circle (2pt); \filldraw (8,4) circle (2pt); \filldraw (9,4) circle (2pt);
\filldraw (1,5) circle (2pt); \filldraw (2,5) circle (2pt); \filldraw (3,5) circle (2pt); \filldraw (4,5) circle (2pt);
\filldraw (5,5) circle (2pt); \filldraw (6,5) circle (2pt); \filldraw (7,5) circle (2pt); \filldraw (8,5) circle (2pt); \filldraw (9,5) circle (2pt); 
\filldraw (1,6) circle (2pt); \filldraw (2,6) circle (2pt); \filldraw (3,6) circle (2pt); \filldraw (4,6) circle (2pt);
\filldraw (5,6) circle (2pt); \filldraw (6,6) circle (2pt); \filldraw (7,6) circle (2pt); \filldraw (8,6) circle (2pt); \filldraw (9,6) circle (2pt);
\filldraw (1,7) circle (2pt); \filldraw (2,7) circle (2pt); \filldraw (3,7) circle (2pt); \filldraw (4,7) circle (2pt);
\filldraw (5,7) circle (2pt); \filldraw (6,7) circle (2pt); \filldraw (7,7) circle (2pt); \filldraw (8,7) circle (2pt); \filldraw (9,7) circle (2pt);

\draw  (1,1)--(7,7); \draw  (2,1)--(8,7); \draw  (3,1)--(9,7); \draw  (4,1)--(9,6); \draw  (5,1)--(9,5); \draw  (6,1)--(9,4); \draw  (7,1)--(9,3); \draw  (8,1)--(9,2);

\draw  (1,2)--(6,7); \draw  (1,3)--(5,7); \draw  (1,4)--(4,7); \draw  (1,5)--(3,7); \draw  (1,6)--(2,7); 

\draw  (1,2)--(2,1); \draw  (1,3)--(3,1); \draw  (1,4)--(4,1); \draw  (1,5)--(5,1); \draw  (1,6)--(6,1); \draw  (1,7)--(7,1);

\draw  (2,7)--(8,1); \draw  (3,7)--(9,1); \draw  (4,7)--(9,2); \draw  (5,7)--(9,3); \draw  (6,7)--(9,4); \draw  (7,7)--(9,5); \draw  (8,7)--(9,6);

\path node at (1.4,1) {$0$}; \path node at (2.4,1) {$3$}; \path node at (3.4,1) {$6$}; \path node at (4.4,1) {$2$}; \path node at (5.4,1) {$5$};
\path node at (6.4,1) {$1$}; \path node at (7.4,1) {$4$}; \path node at (8.4,1) {$0$}; \path node at (9.4,1) {$3$};

\path node at (1.4,2) {$1$}; \path node at (2.4,2) {$4$}; \path node at (3.4,2) {$0$}; \path node at (4.4,2) {$3$}; \path node at (5.4,2) {$6$};
\path node at (6.4,2) {$2$}; \path node at (7.4,2) {$5$}; \path node at (8.4,2) {$1$}; \path node at (9.4,2) {$4$};

\path node at (1.4,3) {$2$}; \path node at (2.4,3) {$5$}; \path node at (3.4,3) {$1$}; \path node at (4.4,3) {$4$}; \path node at (5.4,3) {$0$};
\path node at (6.4,3) {$3$}; \path node at (7.4,3) {$6$}; \path node at (8.4,3) {$2$}; \path node at (9.4,3) {$5$};

\path node at (1.4,4) {$3$}; \path node at (2.4,4) {$6$}; \path node at (3.4,4) {$2$}; \path node at (4.4,4) {$5$}; \path node at (5.4,4) {$1$};
\path node at (6.4,4) {$4$}; \path node at (7.4,4) {$0$}; \path node at (8.4,4) {$3$}; \path node at (9.4,4) {$6$};

\path node at (1.4,5) {$4$}; \path node at (2.4,5) {$0$}; \path node at (3.4,5) {$3$}; \path node at (4.4,5) {$6$}; \path node at (5.4,5) {$2$};
\path node at (6.4,5) {$5$}; \path node at (7.4,5) {$1$}; \path node at (8.4,5) {$4$}; \path node at (9.4,5) {$0$};
 
\path node at (1.4,6) {$5$}; \path node at (2.4,6) {$1$}; \path node at (3.4,6) {$4$}; \path node at (4.4,6) {$0$}; \path node at (5.4,6) {$3$};
\path node at (6.4,6) {$6$}; \path node at (7.4,6) {$2$}; \path node at (8.4,6) {$5$}; \path node at (9.4,6) {$1$};

\path node at (1.4,7) {$6$}; \path node at (2.4,7) {$2$}; \path node at (3.4,7) {$5$}; \path node at (4.4,7) {$1$}; \path node at (5.4,7) {$4$};
\path node at (6.4,7) {$0$}; \path node at (7.4,7) {$3$}; \path node at (8.4,7) {$6$}; \path node at (9.4,7) {$2$};

\path node at (1,0) {$0$}; \path node at (2,0) {$1$}; \path node at (3,0) {$2$}; \path node at (4,0) {$3$}; \path node at (5,0) {$4$};
\path node at (6,0) {$5$}; \path node at (7,0) {$6$}; \path node at (8,0) {$7$};  \path node at (9,0) {$8$};

\path node at (0,1) {$0$}; \path node at (0,2) {$1$}; \path node at (0,3) {$2$}; \path node at (0,4) {$3$}; \path node at (0,5) {$4$};
\path node at (0,6) {$5$}; \path node at (0,7) {$6$};

\path node at (5,-1) {$c)$};


\filldraw (12,1) circle (2pt); \filldraw (13,1) circle (2pt); \filldraw (14,1) circle (2pt); \filldraw (15,1) circle (2pt);
\filldraw (16,1) circle (2pt); \filldraw (17,1) circle (2pt); \filldraw (18,1) circle (2pt); \filldraw (19,1) circle (2pt); \filldraw (20,1) circle (2pt);
\filldraw (12,2) circle (2pt); \filldraw (13,2) circle (2pt); \filldraw (14,2) circle (2pt); \filldraw (15,2) circle (2pt);
\filldraw (16,2) circle (2pt); \filldraw (17,2) circle (2pt); \filldraw (18,2) circle (2pt); \filldraw (19,2) circle (2pt); \filldraw (20,2) circle (2pt);
\filldraw (12,3) circle (2pt); \filldraw (13,3) circle (2pt); \filldraw (14,3) circle (2pt); \filldraw (15,3) circle (2pt); 
\filldraw (16,3) circle (2pt); \filldraw (17,3) circle (2pt); \filldraw (18,3) circle (2pt); \filldraw (19,3) circle (2pt); \filldraw (20,3) circle (2pt);
\filldraw (12,4) circle (2pt); \filldraw (13,4) circle (2pt); \filldraw (14,4) circle (2pt); \filldraw (15,4) circle (2pt);
\filldraw (16,4) circle (2pt); \filldraw (17,4) circle (2pt); \filldraw (18,4) circle (2pt); \filldraw (19,4) circle (2pt); \filldraw (20,4) circle (2pt);
\filldraw (12,5) circle (2pt); \filldraw (13,5) circle (2pt); \filldraw (14,5) circle (2pt); \filldraw (15,5) circle (2pt);
\filldraw (16,5) circle (2pt); \filldraw (17,5) circle (2pt); \filldraw (18,5) circle (2pt); \filldraw (19,5) circle (2pt); \filldraw (20,5) circle (2pt); 
\filldraw (12,6) circle (2pt); \filldraw (13,6) circle (2pt); \filldraw (14,6) circle (2pt); \filldraw (15,6) circle (2pt);
\filldraw (16,6) circle (2pt); \filldraw (17,6) circle (2pt); \filldraw (18,6) circle (2pt); \filldraw (19,6) circle (2pt); \filldraw (20,6) circle (2pt);
\filldraw (12,7) circle (2pt); \filldraw (13,7) circle (2pt); \filldraw (14,7) circle (2pt); \filldraw (15,7) circle (2pt);
\filldraw (16,7) circle (2pt); \filldraw (17,7) circle (2pt); \filldraw (18,7) circle (2pt); \filldraw (19,7) circle (2pt); \filldraw (20,7) circle (2pt);

\draw  (12,1)--(18,7); \draw  (13,1)--(19,7); \draw  (14,1)--(20,7); \draw  (15,1)--(20,6); \draw  (16,1)--(20,5); \draw  (17,1)--(20,4); \draw  (18,1)--(20,3); \draw  (19,1)--(20,2);

\draw  (12,2)--(17,7); \draw  (12,3)--(16,7); \draw  (12,4)--(15,7); \draw  (12,5)--(14,7); \draw  (12,6)--(13,7); 

\draw  (12,2)--(13,1); \draw  (12,3)--(14,1); \draw  (12,4)--(15,1); \draw  (12,5)--(16,1); \draw  (12,6)--(17,1); \draw  (12,7)--(18,1);

\draw  (13,7)--(19,1); \draw  (14,7)--(20,1); \draw  (15,7)--(20,2); \draw  (16,7)--(20,3); \draw  (17,7)--(20,4); \draw  (18,7)--(20,5); \draw  (19,7)--(20,6);

\path node at (12.4,1) {$0$}; \path node at (13.4,1) {$4$}; \path node at (14.4,1) {$1$}; \path node at (15.4,1) {$5$}; \path node at (16.4,1) {$2$};
\path node at (17.4,1) {$6$}; \path node at (18.4,1) {$3$}; \path node at (19.4,1) {$0$}; \path node at (20.4,1) {$4$};

\path node at (12.4,2) {$1$}; \path node at (13.4,2) {$5$}; \path node at (14.4,2) {$2$}; \path node at (15.4,2) {$6$}; \path node at (16.4,2) {$3$};
\path node at (17.4,2) {$0$}; \path node at (18.4,2) {$4$}; \path node at (19.4,2) {$1$}; \path node at (20.4,2) {$5$};

\path node at (12.4,3) {$2$}; \path node at (13.4,3) {$6$}; \path node at (14.4,3) {$3$}; \path node at (15.4,3) {$0$}; \path node at (16.4,3) {$4$};
\path node at (17.4,3) {$1$}; \path node at (18.4,3) {$5$}; \path node at (19.4,3) {$2$}; \path node at (20.4,3) {$6$};

\path node at (12.4,4) {$3$}; \path node at (13.4,4) {$0$}; \path node at (14.4,4) {$4$}; \path node at (15.4,4) {$1$}; \path node at (16.4,4) {$5$};
\path node at (17.4,4) {$2$}; \path node at (18.4,4) {$6$}; \path node at (19.4,4) {$3$}; \path node at (20.4,4) {$0$};

\path node at (12.4,5) {$4$}; \path node at (13.4,5) {$1$}; \path node at (14.4,5) {$5$}; \path node at (15.4,5) {$2$}; \path node at (16.4,5) {$6$};
\path node at (17.4,5) {$3$}; \path node at (18.4,5) {$0$}; \path node at (19.4,5) {$4$}; \path node at (20.4,5) {$1$};
 
\path node at (12.4,6) {$5$}; \path node at (13.4,6) {$2$}; \path node at (14.4,6) {$6$}; \path node at (15.4,6) {$3$}; \path node at (16.4,6) {$0$};
\path node at (17.4,6) {$4$}; \path node at (18.4,6) {$1$}; \path node at (19.4,6) {$5$}; \path node at (20.4,6) {$2$};

\path node at (12.4,7) {$6$}; \path node at (13.4,7) {$3$}; \path node at (14.4,7) {$0$}; \path node at (15.4,7) {$4$}; \path node at (16.4,7) {$1$};
\path node at (17.4,7) {$5$}; \path node at (18.4,7) {$2$}; \path node at (19.4,7) {$6$}; \path node at (20.4,7) {$2$};

\path node at (12,0) {$0$}; \path node at (13,0) {$1$}; \path node at (14,0) {$2$}; \path node at (15,0) {$3$}; \path node at (16,0) {$4$};
\path node at (17,0) {$5$}; \path node at (18,0) {$6$}; \path node at (19,0) {$7$};  \path node at (20,0) {$8$};

\path node at (11,1) {$0$}; \path node at (11,2) {$1$}; \path node at (11,3) {$2$}; \path node at (11,4) {$3$}; \path node at (11,5) {$4$};
\path node at (11,6) {$5$}; \path node at (11,7) {$6$};

\path node at (16,-1) {$d)$};

\filldraw (1,10) circle (2pt); \filldraw (2,10) circle (2pt); \filldraw (3,10) circle (2pt); \filldraw (4,10) circle (2pt);
\filldraw (5,10) circle (2pt); \filldraw (6,10) circle (2pt); \filldraw (7,10) circle (2pt); \filldraw (8,10) circle (2pt); \filldraw (9,10) circle (2pt);
\filldraw (1,11) circle (2pt); \filldraw (2,11) circle (2pt); \filldraw (3,11) circle (2pt); \filldraw (4,11) circle (2pt);
\filldraw (5,11) circle (2pt); \filldraw (6,11) circle (2pt); \filldraw (7,11) circle (2pt); \filldraw (8,11) circle (2pt); \filldraw (9,11) circle (2pt);
\filldraw (1,12) circle (2pt); \filldraw (2,12) circle (2pt); \filldraw (3,12) circle (2pt); \filldraw (4,12) circle (2pt); 
\filldraw (5,12) circle (2pt); \filldraw (6,12) circle (2pt); \filldraw (7,12) circle (2pt); \filldraw (8,12) circle (2pt); \filldraw (9,12) circle (2pt);
\filldraw (1,13) circle (2pt); \filldraw (2,13) circle (2pt); \filldraw (3,13) circle (2pt); \filldraw (4,13) circle (2pt);
\filldraw (5,13) circle (2pt); \filldraw (6,13) circle (2pt); \filldraw (7,13) circle (2pt); \filldraw (8,13) circle (2pt); \filldraw (9,13) circle (2pt);
\filldraw (1,14) circle (2pt); \filldraw (2,14) circle (2pt); \filldraw (3,14) circle (2pt); \filldraw (4,14) circle (2pt);
\filldraw (5,14) circle (2pt); \filldraw (6,14) circle (2pt); \filldraw (7,14) circle (2pt); \filldraw (8,14) circle (2pt); \filldraw (9,14) circle (2pt); 
\filldraw (1,15) circle (2pt); \filldraw (2,15) circle (2pt); \filldraw (3,15) circle (2pt); \filldraw (4,15) circle (2pt);
\filldraw (5,15) circle (2pt); \filldraw (6,15) circle (2pt); \filldraw (7,15) circle (2pt); \filldraw (8,15) circle (2pt); \filldraw (9,15) circle (2pt);
\filldraw (1,16) circle (2pt); \filldraw (2,16) circle (2pt); \filldraw (3,16) circle (2pt); \filldraw (4,16) circle (2pt);
\filldraw (5,16) circle (2pt); \filldraw (6,16) circle (2pt); \filldraw (7,16) circle (2pt); \filldraw (8,16) circle (2pt); \filldraw (9,16) circle (2pt);

 \draw  (1,10)--(7,16); \draw  (2,10)--(8,16); \draw  (3,10)--(9,16); \draw  (4,10)--(9,15); \draw  (5,10)--(9,14); \draw  (6,10)--(9,13); \draw  (7,10)--(9,12); \draw  (8,10)--(9,11);

\draw  (1,11)--(6,16); \draw  (1,12)--(5,16); \draw  (1,13)--(4,16); \draw  (1,14)--(3,16); \draw  (1,15)--(2,16); 

\draw  (1,11)--(2,10); \draw  (1,12)--(3,10); \draw  (1,13)--(4,10); \draw  (1,14)--(5,10); \draw  (1,15)--(6,10); \draw  (1,16)--(7,10);

\draw  (2,16)--(8,10); \draw  (3,16)--(9,10); \draw  (4,16)--(9,11); \draw  (5,16)--(9,12); \draw  (6,16)--(9,13); \draw  (7,16)--(9,14); \draw  (8,16)--(9,15);

\path node at (1.4,10) {$0$}; \path node at (2.4,10) {$1$}; \path node at (3.4,10) {$2$}; \path node at (4.4,10) {$3$}; \path node at (5.4,10) {$4$};
\path node at (6.4,10) {$5$}; \path node at (7.4,10) {$6$}; \path node at (8.4,10) {$0$}; \path node at (9.4,10) {$1$};

\path node at (1.4,11) {$3$}; \path node at (2.4,11) {$4$}; \path node at (3.4,11) {$5$}; \path node at (4.4,11) {$6$}; \path node at (5.4,11) {$0$};
\path node at (6.4,11) {$1$}; \path node at (7.4,11) {$2$}; \path node at (8.4,11) {$3$}; \path node at (9.4,11) {$4$};

\path node at (1.4,12) {$6$}; \path node at (2.4,12) {$0$}; \path node at (3.4,12) {$1$}; \path node at (4.4,12) {$2$}; \path node at (5.4,12) {$3$};
\path node at (6.4,12) {$4$}; \path node at (7.4,12) {$5$}; \path node at (8.4,12) {$6$}; \path node at (9.4,12) {$0$};

\path node at (1.4,13) {$2$}; \path node at (2.4,13) {$3$}; \path node at (3.4,13) {$4$}; \path node at (4.4,13) {$5$}; \path node at (5.4,13) {$6$};
\path node at (6.4,13) {$0$}; \path node at (7.4,13) {$1$}; \path node at (8.4,13) {$2$}; \path node at (9.4,13) {$3$};

\path node at (1.4,14) {$5$}; \path node at (2.4,14) {$6$}; \path node at (3.4,14) {$0$}; \path node at (4.4,14) {$1$}; \path node at (5.4,14) {$2$};
\path node at (6.4,14) {$3$}; \path node at (7.4,14) {$4$}; \path node at (8.4,14) {$5$}; \path node at (9.4,14) {$6$};
 
\path node at (1.4,15) {$1$}; \path node at (2.4,15) {$2$}; \path node at (3.4,15) {$3$}; \path node at (4.4,15) {$4$}; \path node at (5.4,15) {$5$};
\path node at (6.4,15) {$6$}; \path node at (7.4,15) {$0$}; \path node at (8.4,15) {$1$}; \path node at (9.4,15) {$2$};

\path node at (1.4,16) {$4$}; \path node at (2.4,16) {$5$}; \path node at (3.4,16) {$6$}; \path node at (4.4,16) {$0$}; \path node at (5.4,16) {$1$};
\path node at (6.4,16) {$2$}; \path node at (7.4,16) {$3$}; \path node at (8.4,16) {$4$}; \path node at (9.4,16) {$5$};

\path node at (1,9) {$0$}; \path node at (2,9) {$1$}; \path node at (3,9) {$2$}; \path node at (4,9) {$3$}; \path node at (5,9) {$4$};
\path node at (6,9) {$5$}; \path node at (7,9) {$6$}; \path node at (8,9) {$7$};  \path node at (9,9) {$8$};

\path node at (0,10) {$0$}; \path node at (0,11) {$1$}; \path node at (0,12) {$2$}; \path node at (0,13) {$3$}; \path node at (0,14) {$4$};
\path node at (0,15) {$5$}; \path node at (0,16) {$6$};

\path node at (5,8) {$a)$};

\filldraw (12,10) circle (2pt); \filldraw (13,10) circle (2pt); \filldraw (14,10) circle (2pt); \filldraw (15,10) circle (2pt);
\filldraw (16,10) circle (2pt); \filldraw (17,10) circle (2pt); \filldraw (18,10) circle (2pt); \filldraw (19,10) circle (2pt); \filldraw (20,10) circle (2pt);
\filldraw (12,11) circle (2pt); \filldraw (13,11) circle (2pt); \filldraw (14,11) circle (2pt); \filldraw (15,11) circle (2pt);
\filldraw (16,11) circle (2pt); \filldraw (17,11) circle (2pt); \filldraw (18,11) circle (2pt); \filldraw (19,11) circle (2pt); \filldraw (20,11) circle (2pt);
\filldraw (12,12) circle (2pt); \filldraw (13,12) circle (2pt); \filldraw (14,12) circle (2pt); \filldraw (15,12) circle (2pt); 
\filldraw (16,12) circle (2pt); \filldraw (17,12) circle (2pt); \filldraw (18,12) circle (2pt); \filldraw (19,12) circle (2pt); \filldraw (20,12) circle (2pt);
\filldraw (12,13) circle (2pt); \filldraw (13,13) circle (2pt); \filldraw (14,13) circle (2pt); \filldraw (15,13) circle (2pt);
\filldraw (16,13) circle (2pt); \filldraw (17,13) circle (2pt); \filldraw (18,13) circle (2pt); \filldraw (19,13) circle (2pt); \filldraw (20,13) circle (2pt);
\filldraw (12,14) circle (2pt); \filldraw (13,14) circle (2pt); \filldraw (14,14) circle (2pt); \filldraw (15,14) circle (2pt);
\filldraw (16,14) circle (2pt); \filldraw (17,14) circle (2pt); \filldraw (18,14) circle (2pt); \filldraw (19,14) circle (2pt); \filldraw (20,14) circle (2pt); 
\filldraw (12,15) circle (2pt); \filldraw (13,15) circle (2pt); \filldraw (14,15) circle (2pt); \filldraw (15,15) circle (2pt);
\filldraw (16,15) circle (2pt); \filldraw (17,15) circle (2pt); \filldraw (18,15) circle (2pt); \filldraw (19,15) circle (2pt); \filldraw (20,15) circle (2pt);
\filldraw (12,16) circle (2pt); \filldraw (13,16) circle (2pt); \filldraw (14,16) circle (2pt); \filldraw (15,16) circle (2pt);
\filldraw (16,16) circle (2pt); \filldraw (17,16) circle (2pt); \filldraw (18,16) circle (2pt); \filldraw (19,16) circle (2pt); \filldraw (20,16) circle (2pt);

\draw  (12,10)--(18,16); \draw  (13,10)--(19,16); \draw  (14,10)--(20,16); \draw  (15,10)--(20,15); \draw  (16,10)--(20,14); \draw  (17,10)--(20,13); \draw  (18,10)--(20,12); \draw  (19,10)--(20,11);

\draw  (12,11)--(17,16); \draw  (12,12)--(16,16); \draw  (12,13)--(15,16); \draw  (12,14)--(14,16); \draw  (12,15)--(13,16); 

\draw  (12,11)--(13,10); \draw  (12,12)--(14,10); \draw  (12,13)--(15,10); \draw  (12,14)--(16,10); \draw  (12,15)--(17,10); \draw  (12,16)--(18,10);

\draw  (13,16)--(19,10); \draw  (14,16)--(20,10); \draw  (15,16)--(20,11); \draw  (16,16)--(20,12); \draw  (17,16)--(20,13); \draw  (18,16)--(20,14); \draw  (19,16)--(20,15);

\path node at (12.4,10) {$0$}; \path node at (13.4,10) {$1$}; \path node at (14.4,10) {$2$}; \path node at (15.4,10) {$3$}; \path node at (16.4,10) {$4$};
\path node at (17.4,10) {$5$}; \path node at (18.4,10) {$6$}; \path node at (19.4,10) {$0$}; \path node at (20.4,10) {$1$};

\path node at (12.4,11) {$4$}; \path node at (13.4,11) {$5$}; \path node at (14.4,11) {$6$}; \path node at (15.4,11) {$0$}; \path node at (16.4,11) {$1$};
\path node at (17.4,11) {$2$}; \path node at (18.4,11) {$3$}; \path node at (19.4,11) {$4$}; \path node at (20.4,11) {$5$};

\path node at (12.4,12) {$1$}; \path node at (13.4,12) {$2$}; \path node at (14.4,12) {$3$}; \path node at (15.4,12) {$4$}; \path node at (16.4,12) {$5$};
\path node at (17.4,12) {$6$}; \path node at (18.4,12) {$0$}; \path node at (19.4,12) {$1$}; \path node at (20.4,12) {$2$};

\path node at (12.4,13) {$5$}; \path node at (13.4,13) {$6$}; \path node at (14.4,13) {$0$}; \path node at (15.4,13) {$1$}; \path node at (16.4,13) {$2$};
\path node at (17.4,13) {$3$}; \path node at (18.4,13) {$4$}; \path node at (19.4,13) {$5$}; \path node at (20.4,13) {$6$};

\path node at (12.4,14) {$2$}; \path node at (13.4,14) {$3$}; \path node at (14.4,14) {$4$}; \path node at (15.4,14) {$5$}; \path node at (16.4,14) {$6$};
\path node at (17.4,14) {$0$}; \path node at (18.4,14) {$1$}; \path node at (19.4,14) {$2$}; \path node at (20.4,14) {$3$};
 
\path node at (12.4,15) {$6$}; \path node at (13.4,15) {$0$}; \path node at (14.4,15) {$1$}; \path node at (15.4,15) {$2$}; \path node at (16.4,15) {$3$};
\path node at (17.4,15) {$4$}; \path node at (18.4,15) {$5$}; \path node at (19.4,15) {$6$}; \path node at (20.4,15) {$0$};

\path node at (12.4,16) {$3$}; \path node at (13.4,16) {$4$}; \path node at (14.4,16) {$5$}; \path node at (15.4,16) {$6$}; \path node at (16.4,16) {$0$};
\path node at (17.4,16) {$1$}; \path node at (18.4,16) {$2$}; \path node at (19.4,16) {$3$}; \path node at (20.4,16) {$4$};

\path node at (12,9) {$0$}; \path node at (13,9) {$1$}; \path node at (14,9) {$2$}; \path node at (15,9) {$3$}; \path node at (16,9) {$4$};
\path node at (17,9) {$5$}; \path node at (18,9) {$6$}; \path node at (19,9) {$7$};  \path node at (20,9) {$8$};

\path node at (11,10) {$0$}; \path node at (11,11) {$1$}; \path node at (11,12) {$2$}; \path node at (11,13) {$3$}; \path node at (11,14) {$4$};
\path node at (11,15) {$5$}; \path node at (11,16) {$6$};

\path node at (16,8) {$b)$};	}	
\end{tikzpicture}
\caption{Four  $L(2,1)$-labelings of $P_9\times P_7$ that determined $L(2,1)$-labelings of direct graph bundle $C_9\times^{\sigma_\ell}C_7$ according to the cyclic $\ell$-shift $\sigma_\ell$: \\
a) $ \ell =3, f_1(i,j)=\left[i+3j\right]\mod 7$,\\
b) $ \ell =4, f_2(i,j)=\left[i+4j)\right]\mod 7 $, \\
c) $ \ell=6, g_1(i,j)=\left[3i+j\right]\mod 7 $, \\
d) $ \ell=1, g_2(i,j)=\left[4i+j\right]\mod 7$.}
\label{primer1}
\end{center}
\end{figure}

\newpage

\section{Cartesian graph bundle $C_m\Box^{\sigma_\ell} C_n$  }


\begin{theorem}
\label{izrek2}
Let $d\geq 1$, $m\geq 3$ and $s=2d+3$. Let $X= C_m\Box^{\sigma_\ell} C_n$ be a Cartesian graph bundle with fibre $C_n$ and base $C_m$ and let $n$ be a multiple of $s$.  Then $\lambda^d_1(X)\leq 2d+2$, with equality if $1\leq d\leq 4$, if one of the following statements holds:
\begin{description}
\item[(a)]  $\ell$ has a form of $\ell =[ks+(-1)^a2dm]\mod n$, where $a\in\{1,2\}$ and $k\in \ZZ.$  
\item[(b)]  $d=3t+2$ for some $t\in \{0,1,\ldots\}$ and $\ell$ has a form of $$\ell =\left[ks-(2t+3)(d+a)m\right]\mod n, {\mathop{\rm where\;}} a\in\{1,2\}\;{\mathop{\rm and\;}}k\in \ZZ.$$  
\item[(c)]  $d=3t+1$ for some $t\in \{0,1,\ldots\}$ and $\ell$ has a form of $$\ell =\left[ks+(2t+1)(d+a)m\right]\mod n,{\mathop{\rm where\;}} a\in\{1,2\}\;{\mathop{\rm and\;}}k\in \ZZ.$$    
\item[(d)]  $d=3t$ for some $t\in\{1,2,\ldots\}$,  $m=ps+3t'\geq 3$ for some $p\in \{0,1,\ldots\}$ and some $t'\in\{0,1,\ldots ,2t\}$ and  $\ell$ has a form of $$\ell =\left[ks+(i+a-1)\frac{s}{3}-(-1)^at'\right]\mod n, {\mathop{\rm where\;}} a\in\{1,2\},\; k\in \ZZ\; {\mathop{\rm and\;}}  i\in\{0,1,2\}.$$
\end{description} 
\end{theorem}

\noindent {\bf Proof:}

We will present four  $L(d,1)$-labelings of $C_m\Box^{\sigma_\ell} C_n$ using labels $0,1,\ldots , 2d+2$ according to the cyclic $\ell$-shift $\sigma_\ell$. Let $v=(i,j)\in V(X)$. 
\begin{enumerate} 
\item  Assume that the statement (a) from the theorem holds. 
For $a\in\{1,2\}$ define labeling $f_a$ of $v$ as  
$$f_a(v)=\left[di+(d+a)j\right]\mod s.$$ 

\item Assume that one of the statements (b), (c) or (d) holds. For $a\in\{1,2\}$ define labeling $g_a$ of $v$ as  
$$g_a(v)=\left[(d+a)i+dj\right]\mod s.$$ 

\end{enumerate}

All assignments are clearly well-defined. 


Since $ C_m\Box^{\sigma_\ell} C_n$ can be represented  as the graph obtained from the product $P_m\Box C_n$ by adding edges between vertex sets $\{m-1\}\times V(C_n)$ and $\{0\}\times V(C_n)$ first observe the product $P_m\Box C_n$. 

Let $X$ have labeling $f_a$. Then $v$ be assigned the integer  $$f_a(v)=\left[di+(d+a)j\right]\mod s.$$ 

Let $w$ be a vertex adjacent to $v$ in $P_m\Box C_n$. Then $v$ and $w$ differ in exactly one coordinate. More precisely, $w$ is of the form $w= (i+i',(j+j')\mod n)$ with one of the following properties:
\begin{enumerate}
\item if $i=0$, then $i'=1,j'=0$  or $i'=0,|j'|= 1$, 
\item if $i=m-1$, then $i'=-1,j'=0$  or $|j'|= 1, i'=0$, 
\item if $i\in \{1,2,\ldots m-2\}$, then $|i'|=1,j'=0$  or $|j'|= 1, i'=0$.
\end{enumerate}

It is clear that $$f_a(w)=\left[di+(d+a)j+di'+(d+a)j'\right]\mod s.$$

By Corollary \ref{posledica 1}, to show that $|f_a(v)-f_a(w)|\geq d$, it is enough to show that $$d\leq \left[|di'+(d+a)j'|\mod s\right] \leq s-d=d+3$$ 
and, by symmetry, if $X$ has labeling $g_a$, that $$d\leq \left[|(d+a)i'+dj'|\mod s\right] \leq d+3.$$ Both claims are true since $|di'+(d+a)j'|\mod s,|(d+a)i'+j'|\mod s\in \{d,d+1,d+2,d+3\}$.

Next, let $z$ be a vertex at a distance of two from $v$ in $P_m\Box C_n$. Then $z$ is of the form $z= (i+i'',(j+j'')\mod n)$, where $v$ and $z$ differ in exactly one (a) or both coordinates (b). In all cases there are three options for the values of $i''$ and $j''$. In (a):
\begin{enumerate}
\item if  $i\in \{0,1\}$, then  $i''=0,|j''|=2$ or $i''=2$ and $j''=0$,
\item if $i\in \{m-2,m-1\}$, then $i''=0,|j''|=2$ or $i''=-2$ and $j''=0$,
\item otherwise, $i''=0,|j''|=2$ or $|i''|=2,j''=0 $.
\end{enumerate}

And in (b):
\begin{enumerate}
\item if  $i=0$, then  $i''=1,|j''|=1$,
\item if $i=m-1$, then $i''=-1,|j''|=1$,
\item otherwise, $|i''|=1,|j''|=1$.
\end{enumerate}

Let $X$ have labeling $f_a$. Then $$f_a(z)=\left[di+(d+a)j+di''+(d+a)j''\right]\mod s.$$ Since we want to prove that $|f_a(v)- f_a(z)|\geq 1$, it is enough to prove, by Corollary \ref{posledica 1},  that  $$1\leq [|di''+(d+a)j''|\mod s]\leq 2d+2$$ and if $X$ has the labeling $g_a$, that $$1\leq [|(d+a)i''+dj''|\mod s]\leq 2d+2.$$  
The reader can verify that in (a) we get  $|di''+(d+a)j''|\mod s,|(d+a)i''+dj''|\mod s\in\{1,3, 2d,2d+2\}$ and in (b) $|di''+(d+a)j''|\mod s, |(d+a)i''+dj''|\mod s\in \{ 1,2,2d+1,2d+2\}$. Hence, both results follow.

In the following, we are interested in edges  in $X= C_m\Box^{\sigma_\ell} C_n$ between the fiber over vertex $m-1$ and the fiber over vertex $0$. These are edges $(m-1,u)(0, (u+\ell)\mod n), u\in V(C_n)$ or in another form $(0,u)(m-1,(u-\ell)\mod n), u\in V(C_n)$. 

Let first observe two adjacent vertices  $v$ and $w$. Let $X$ have labeling $f_a$ and let $v=(m-1,j)$. Then $w$ is of the form  $w=(0,(j+\ell)\mod n)$, where $\ell =[ks+(-1)^a 2dm]\mod n$ for some $k\in \ZZ$. Vertex $v$ be assigned the integer  $$f_a(v)=\left[d(m-1)+(d+a)j\right]\mod s$$ and $w$ be assigned the integer $$f_a(w)=\left[(d+a)(j+\ell)\right]\mod s.$$  

We claim that $|f_a(v)-f_a(w)|\geq d$ or, by Corollary \ref{posledica 1},  that 
\begin{equation}
\label{c}
d\leq [|d(m-1)-(d+a)\ell|\mod s]\leq d+3.
\end{equation}


Since there exists $k'\in \ZZ$ such that $\ell =k's+(-1)^a 2dm$ it follows that $$|d(m-1)-(d+a)\ell|\mod s=|d(m-1)-(d+a)(k's+(-1)^a 2dm)|\mod s=$$ $$=|dm(1-(-1)^a2d-(-1)^a2a)-d-(d+a)k's|\mod s=$$  $$=|dm(-1)^{a+1}s-d-(d+a)k's|\mod s=|s(dm(-1)^{a+1}-(d+a)k')-d|\mod s $$ (since $1-(-1)^a2d-(-1)^a2a$ equals $s$ for $a=1$ and  $-s$ for $a=2$). By Corollary \ref{posledica 2},  this is equal to $d$ or to $s-d=d+3$ and the claim is true. 

Further, observe labeling $g_a$. Then one of the statements (b), (c) or (d) of the Theorem applies. Let these be Cases b), c) and d), accordingly. Since $$g_a(v)=\left[(d+a)(m-1)+dj\right]\mod s$$ and $$g_a(w)=\left[d(j+\ell)\right]\mod s, $$ 
in all three cases we have to show that 
\begin{equation}
\label{d}
d\leq [|(d+a)(m-1)-d\ell|\mod s]\leq d+3,
\end{equation} 

 by Corollary \ref{posledica 1}.

{\bf Case b)}
Let $d=3t+2$ (then $s=6t+7$). 
Since $\ell=k's-(2t+3)(d+a)m$  for some $k'\in \ZZ$ we get  $$|(d+a)(m-1)-d\ell|\mod s=|(d+a)(m-1)-d(k's-(2t+3)(d+a)m)|\mod s=$$ $$|(d+a)m(1+d(2t+3))-(d+a)-dk's|\mod s=$$ $$|(d+a)m(t+1)s-(d+a)-dk's|\mod s=|s((d+a)m(t+1)-dk')-(d+a)|\mod s.$$ 
By Corollary \ref{posledica 2}, this is equal to $(d+a)\mod s, s-((d+a)\mod s)\in \{d+1,d+2\}$.

{\bf Case c)}
Let $d=3t+1$ (then $s=6t+5$).
Since $\ell=k's+(2t+1)(d+a)m$  for some $k'\in \ZZ$ it holds $$|(d+a)(m-1)-d\ell|\mod s=|(d+a)(m-1)-
d(k's+(2t+1)(d+a)m)|\mod s=$$ $$|(d+a)m(1-d(2t+1))-(d+a)-dk's|\mod s=$$ $$|(d+a)m(-ts)-(d+a)-dk's|\mod s=|s(-mt(d+a)-dk')-(d+a)|\mod s.$$
As in case a) this is equal to $d+1$ or $d+2$.

{\bf Case d)}
Let $d=3t$ (then $s=6t+3$). 
Since $\ell =k's+(i+a-1)\frac{s}{3}-(-1)^at'$ for some $k'\in \ZZ$ it holds
$$|(d+a)(m-1)-d\ell|\mod s=$$ $$|(d+a)(ps+3t'-1)-d(k's+(i+a-1)\frac{s}{3}-(-1)^at')|\mod s =$$   $$|s[p(d+a)-dk'-t(i+a-1)]+(d+a)(3t'-1)+d(-1)^at'|\mod s =$$ 
$$|s[p(d+a)-dk'-t(i+a-1)]+t'(3(d+a)+d(-1)^a)-(d+a)|\mod s.$$
For $a=1$ ve get $$|s[p(d+1)-dk'-ti)]+st'-(d+1)|\mod s$$ and for $a=2$ we get $$|s[p(d+2)-dk'-t(i+1)]+2st'-(d+2)|\mod s .$$ The desired result follows, since in this case we also get $d+1$ or $d+2$.

Finally, observe vertices at a distance of two, where the shortest path between them contains some edge between the fiber over $m-1$ and the fiber over $0$. Let $v$ and $z$ be two such vertices. Then $v$ and $z$ are from adjacent fibers (from fiber over $m-1$ and over $0$) or from non-adjacent fibers (from fiber over $m-2$ and over $0$ or from fiber over $m-1$ and  over $1$). Let us first observe the case  where $v$ and $z$ are from non-adjacent fibers. 
 
\begin{enumerate}
\item Let $v=(m-2,j)$. Then $z$ is of the form $z=(0, (j+\ell)\mod n)$, where $\ell =[ks+(-1)^a2dm]\mod n$ for some $k\in \ZZ$ (similarly, $v=(m-1,j)$ and $z=(1, (j+\ell)\mod n)$). Let $X$ have labeling $f_a$. Vertex $v$ be assigned the integer $$f_a(v)=\left[d(m-2)+(d+a)j)\right]\mod s$$ and vertex $z$  the integer $$f_a(z)=\left[(d+a)(j+\ell)\right]\mod s.$$ 

We claim that $|f_a(v)-f_a(z)|\geq 1$ or, by Corollary \ref{posledica 1}, that  $$1\leq \left[|d(m-2)-(d+a)\ell|\mod s\right] \leq 2d+2.$$

In the proof of (\ref{c}), we showed that $d(m-1)-(d+a)\ell=k''s-d$ for some $k''\in \ZZ$. 
Therefore $$|d(m-2)-(d+a)\ell|\mod s=|d(m-1)-d-(d+a)\ell|\mod s=|k''s-2d|\mod s.$$ By Corollary \ref{posledica 2},  this is equal to $2d$ or to $s-2d=3$ and the desired follows. \\

Let now $X$ have labeling $g_a$.  We need to show that 
$$1\leq [|(d+a)(m-2)-d\ell|\mod s]\leq 2d+2.$$ 

When we proved (\ref{d}), in all three cases we found that $(d+a)(m-1)-d\ell=k''s-(d+a)$ for some $k''\in \ZZ$. Therefore we get $$|(d+a)(m-2)-d\ell|\mod s=|(d+a)(m-1)-(d+a)-d\ell|\mod s=$$ $$|k''s-2(d+a)|\mod s.$$ Since  $(2(d+a))\mod s,s-((2(d+a))\mod s) \in \{1,2d+2\},$ 
the desired result follows.

\item Let $v$ and $z$ be from adjacent fibers and let  $v=(m-1,j)$. Let $X$ have the labeling $f_a$. Then $z$ is of the form $z=(0, (j+j'+\ell)\mod n)$, where $|j'|=1$ and  $\ell =[ks+(-1)^a2dm]\mod n$ for some $k\in \ZZ$. In this case  $$f_a(v)=\left[d(m-1)+(d+a)j)\right]\mod s$$ and $$f_a(z)=\left[(d+a)(j+j'+\ell)\right]\mod s.$$ Again, it is enough to show that
$$1\leq [|d(m-1)-(d+a)(j'+\ell)|\mod s]\leq 2d+2.$$  Recall that $d(m-1)-(d+a)\ell=k''s-d$ for some $k''\in \ZZ$ (see the proof of (\ref{c})). Therefore $$|d(m-1)-(d+a)(j'+\ell)|\mod s=|d(m-1)-(d+a)j'-(d+a)\ell|\mod s=$$ $$|k''s-(d+(d+a)j')|\mod s.$$ The reader can verify that $(d+(d+a)j')\mod s\in \{2d+1,2d+2\}$ and $s-((d+(d+a)j')\mod s)\in \{1,2\}$, so the desired result follows.\\

Let now $X$ have the labeling $g_a$. We need to show that $$1\leq [|(d+a)(m-1)-d(j'+\ell)|\mod s]\leq 2d+2.$$

Recall that  $(d+a)(m-1)-d\ell=k''s-(d+a)$ for some $k''\in \ZZ$ (see the proof of (\ref{d})). Therefore $$|(d+a)(m-1)-d(j'+\ell)|\mod s=|k''s-((d+a)+dj')|\mod s.$$ Since $((d+a)+dj')\mod s, s-(((d+a)+dj')\mod s)\in \{1,2,2d+1,2d+2\}$, the desired result follows. 

Accordingly, two vertices in $C_m\Box^{\sigma_\ell} C_n$ that are at distance of two from each other  receive different labels.

\end{enumerate}

We showed that $\lambda^d_1(C_m\Box^{\sigma_\ell} C_n)\leq 2d+2$. Further, $C_m\Box^{\sigma_\ell} C_n$ being a regular graph of degree 4, an application of Lemma \ref{lemma 1} to the preceding statement shows that $\lambda^d_1(C_m\Box^{\sigma_\ell} C_n)= 2d+2$, if $1\leq d \leq 4$.

\rule{2mm}{2mm}

The foregoing scheme is illustrated in Fig. \ref{slika2}, where an $L(2,1)$-labeling of $P_9\Box P_7$ appears toward that of  $C_9\Box^{\sigma_\ell}C_7$ for $\ell=1,3,4,6$.


\begin{figure}[h]
\label{slika2}
\begin{center}
\begin{tikzpicture} [scale=0.65]

	{    				
\filldraw (1,1) circle (2pt); \filldraw (2,1) circle (2pt); \filldraw (3,1) circle (2pt); \filldraw (4,1) circle (2pt);
\filldraw (5,1) circle (2pt); \filldraw (6,1) circle (2pt); \filldraw (7,1) circle (2pt); \filldraw (8,1) circle (2pt); \filldraw (9,1) circle (2pt);
\filldraw (1,2) circle (2pt); \filldraw (2,2) circle (2pt); \filldraw (3,2) circle (2pt); \filldraw (4,2) circle (2pt);
\filldraw (5,2) circle (2pt); \filldraw (6,2) circle (2pt); \filldraw (7,2) circle (2pt); \filldraw (8,2) circle (2pt); \filldraw (9,2) circle (2pt);
\filldraw (1,3) circle (2pt); \filldraw (2,3) circle (2pt); \filldraw (3,3) circle (2pt); \filldraw (4,3) circle (2pt); 
\filldraw (5,3) circle (2pt); \filldraw (6,3) circle (2pt); \filldraw (7,3) circle (2pt); \filldraw (8,3) circle (2pt); \filldraw (9,3) circle (2pt);
\filldraw (1,4) circle (2pt); \filldraw (2,4) circle (2pt); \filldraw (3,4) circle (2pt); \filldraw (4,4) circle (2pt);
\filldraw (5,4) circle (2pt); \filldraw (6,4) circle (2pt); \filldraw (7,4) circle (2pt); \filldraw (8,4) circle (2pt); \filldraw (9,4) circle (2pt);
\filldraw (1,5) circle (2pt); \filldraw (2,5) circle (2pt); \filldraw (3,5) circle (2pt); \filldraw (4,5) circle (2pt);
\filldraw (5,5) circle (2pt); \filldraw (6,5) circle (2pt); \filldraw (7,5) circle (2pt); \filldraw (8,5) circle (2pt); \filldraw (9,5) circle (2pt); 
\filldraw (1,6) circle (2pt); \filldraw (2,6) circle (2pt); \filldraw (3,6) circle (2pt); \filldraw (4,6) circle (2pt);
\filldraw (5,6) circle (2pt); \filldraw (6,6) circle (2pt); \filldraw (7,6) circle (2pt); \filldraw (8,6) circle (2pt); \filldraw (9,6) circle (2pt);
\filldraw (1,7) circle (2pt); \filldraw (2,7) circle (2pt); \filldraw (3,7) circle (2pt); \filldraw (4,7) circle (2pt);
\filldraw (5,7) circle (2pt); \filldraw (6,7) circle (2pt); \filldraw (7,7) circle (2pt); \filldraw (8,7) circle (2pt); \filldraw (9,7) circle (2pt);

\draw  (1,1)--(9,1); \draw  (1,2)--(9,2);\draw  (1,3)--(9,3);\draw  (1,4)--(9,4);\draw  (1,5)--(9,5);\draw  (1,6)--(9,6);\draw  (1,7)--(9,7);

\draw  (1,1)--(1,7); \draw  (2,1)--(2,7);\draw  (3,1)--(3,7);\draw  (4,1)--(4,7);\draw  (5,1)--(5,7);\draw  (6,1)--(6,7);\draw  (7,1)--(7,7);
\draw  (8,1)--(8,7);\draw  (9,1)--(9,7);

\path node at (1.3,1.3) {$0$}; \path node at (2.3,1.3) {$3$}; \path node at (3.4,1.3) {$6$}; \path node at (4.4,1.3) {$2$}; \path node at (5.4,1.3) {$5$};
\path node at (6.4,1.3) {$1$}; \path node at (7.4,1.3) {$4$}; \path node at (8.4,1.3) {$0$}; \path node at (9.4,1.3) {$3$};

\path node at (1.3,2.3) {$2$}; \path node at (2.3,2.3) {$5$}; \path node at (3.3,2.3) {$1$}; \path node at (4.3,2.3) {$4$}; \path node at (5.3,2.3) {$0$};
\path node at (6.4,2.3) {$3$}; \path node at (7.4,2.3) {$6$}; \path node at (8.4,2.3) {$2$}; \path node at (9.4,2.3) {$5$};

\path node at (1.3,3.3) {$4$}; \path node at (2.3,3.3) {$0$}; \path node at (3.3,3.3) {$3$}; \path node at (4.3,3.3) {$6$}; \path node at (5.3,3.3) {$2$};
\path node at (6.3,3.3) {$5$}; \path node at (7.3,3.3) {$1$}; \path node at (8.3,3.3) {$4$}; \path node at (9.3,3.3) {$0$};

\path node at (1.3,4.3) {$6$}; \path node at (2.3,4.3) {$2$}; \path node at (3.3,4.3) {$5$}; \path node at (4.3,4.3) {$1$}; \path node at (5.3,4.3) {$4$};
\path node at (6.3,4.3) {$0$}; \path node at (7.3,4.3) {$3$}; \path node at (8.3,4.3) {$6$}; \path node at (9.3,4.3) {$2$};

\path node at (1.3,5.3) {$1$}; \path node at (2.3,5.3) {$4$}; \path node at (3.3,5.3) {$0$}; \path node at (4.3,5.3) {$3$}; \path node at (5.3,5.3) {$6$};
\path node at (6.3,5.3) {$2$}; \path node at (7.3,5.3) {$5$}; \path node at (8.3,5.3) {$1$}; \path node at (9.3,5.3) {$4$};
 
\path node at (1.3,6.3) {$3$}; \path node at (2.3,6.3) {$6$}; \path node at (3.3,6.3) {$2$}; \path node at (4.3,6.3) {$5$}; \path node at (5.3,6.3) {$1$};
\path node at (6.3,6.3) {$4$}; \path node at (7.3,6.3) {$0$}; \path node at (8.3,6.3) {$3$}; \path node at (9.3,6.3) {$6$};

\path node at (1.3,7.3) {$5$}; \path node at (2.3,7.3) {$1$}; \path node at (3.3,7.3) {$4$}; \path node at (4.3,7.3) {$0$}; \path node at (5.3,7.3) {$3$};
\path node at (6.3,7.3) {$6$}; \path node at (7.3,7.3) {$2$}; \path node at (8.3,7.3) {$5$}; \path node at (9.3,7.3) {$1$};

\path node at (1,0) {$0$}; \path node at (2,0) {$1$}; \path node at (3,0) {$2$}; \path node at (4,0) {$3$}; \path node at (5,0) {$4$};
\path node at (6,0) {$5$}; \path node at (7,0) {$6$}; \path node at (8,0) {$7$};  \path node at (9,0) {$8$};

\path node at (0,1) {$0$}; \path node at (0,2) {$1$}; \path node at (0,3) {$2$}; \path node at (0,4) {$3$}; \path node at (0,5) {$4$};
\path node at (0,6) {$5$}; \path node at (0,7) {$6$};

\path node at (5,-1) {$c)$};


\filldraw (12,1) circle (2pt); \filldraw (13,1) circle (2pt); \filldraw (14,1) circle (2pt); \filldraw (15,1) circle (2pt);
\filldraw (16,1) circle (2pt); \filldraw (17,1) circle (2pt); \filldraw (18,1) circle (2pt); \filldraw (19,1) circle (2pt); \filldraw (20,1) circle (2pt);
\filldraw (12,2) circle (2pt); \filldraw (13,2) circle (2pt); \filldraw (14,2) circle (2pt); \filldraw (15,2) circle (2pt);
\filldraw (16,2) circle (2pt); \filldraw (17,2) circle (2pt); \filldraw (18,2) circle (2pt); \filldraw (19,2) circle (2pt); \filldraw (20,2) circle (2pt);
\filldraw (12,3) circle (2pt); \filldraw (13,3) circle (2pt); \filldraw (14,3) circle (2pt); \filldraw (15,3) circle (2pt); 
\filldraw (16,3) circle (2pt); \filldraw (17,3) circle (2pt); \filldraw (18,3) circle (2pt); \filldraw (19,3) circle (2pt); \filldraw (20,3) circle (2pt);
\filldraw (12,4) circle (2pt); \filldraw (13,4) circle (2pt); \filldraw (14,4) circle (2pt); \filldraw (15,4) circle (2pt);
\filldraw (16,4) circle (2pt); \filldraw (17,4) circle (2pt); \filldraw (18,4) circle (2pt); \filldraw (19,4) circle (2pt); \filldraw (20,4) circle (2pt);
\filldraw (12,5) circle (2pt); \filldraw (13,5) circle (2pt); \filldraw (14,5) circle (2pt); \filldraw (15,5) circle (2pt);
\filldraw (16,5) circle (2pt); \filldraw (17,5) circle (2pt); \filldraw (18,5) circle (2pt); \filldraw (19,5) circle (2pt); \filldraw (20,5) circle (2pt); 
\filldraw (12,6) circle (2pt); \filldraw (13,6) circle (2pt); \filldraw (14,6) circle (2pt); \filldraw (15,6) circle (2pt);
\filldraw (16,6) circle (2pt); \filldraw (17,6) circle (2pt); \filldraw (18,6) circle (2pt); \filldraw (19,6) circle (2pt); \filldraw (20,6) circle (2pt);
\filldraw (12,7) circle (2pt); \filldraw (13,7) circle (2pt); \filldraw (14,7) circle (2pt); \filldraw (15,7) circle (2pt);
\filldraw (16,7) circle (2pt); \filldraw (17,7) circle (2pt); \filldraw (18,7) circle (2pt); \filldraw (19,7) circle (2pt); \filldraw (20,7) circle (2pt);

\draw  (12,1)--(20,1); \draw  (12,2)--(20,2);\draw  (12,3)--(20,3);\draw  (12,4)--(20,4);\draw  (12,5)--(20,5);\draw  (12,6)--(20,6);\draw  (12,7)--(20,7);

\draw  (12,1)--(12,7); \draw  (13,1)--(13,7);\draw  (14,1)--(14,7);\draw  (15,1)--(15,7);\draw  (16,1)--(16,7);\draw  (17,1)--(17,7);\draw  (18,1)--(18,7);
\draw  (19,1)--(19,7);\draw  (20,1)--(20,7);

\path node at (12.3,1.3) {$0$}; \path node at (13.3,1.3) {$4$}; \path node at (14.3,1.3) {$1$}; \path node at (15.3,1.3) {$5$}; \path node at (16.3,1.3) {$2$};
\path node at (17.3,1.3) {$6$}; \path node at (18.3,1.3) {$3$}; \path node at (19.3,1.3) {$0$}; \path node at (20.3,1.3) {$4$};

\path node at (12.3,2.3) {$2$}; \path node at (13.3,2.3) {$6$}; \path node at (14.3,2.3) {$3$}; \path node at (15.3,2.3) {$0$}; \path node at (16.3,2.3) {$4$};
\path node at (17.3,2.3) {$1$}; \path node at (18.3,2.3) {$5$}; \path node at (19.3,2.3) {$2$}; \path node at (20.3,2.3) {$6$};

\path node at (12.3,3.3) {$4$}; \path node at (13.3,3.3) {$1$}; \path node at (14.3,3.3) {$5$}; \path node at (15.3,3.3) {$2$}; \path node at (16.3,3.3) {$6$};
\path node at (17.3,3.3) {$3$}; \path node at (18.3,3.3) {$0$}; \path node at (19.3,3.3) {$4$}; \path node at (20.3,3.3) {$1$};

\path node at (12.3,4.3) {$6$}; \path node at (13.3,4.3) {$3$}; \path node at (14.3,4.3) {$0$}; \path node at (15.3,4.3) {$4$}; \path node at (16.3,4.3) {$1$};
\path node at (17.3,4.3) {$5$}; \path node at (18.3,4.3) {$2$}; \path node at (19.3,4.3) {$6$}; \path node at (20.3,4.3) {$3$};

\path node at (12.3,5.3) {$1$}; \path node at (13.3,5.3) {$5$}; \path node at (14.3,5.3) {$2$}; \path node at (15.3,5.3) {$6$}; \path node at (16.3,5.3) {$3$};
\path node at (17.3,5.3) {$0$}; \path node at (18.3,5.3) {$4$}; \path node at (19.3,5.3) {$1$}; \path node at (20.3,5.3) {$5$};
 
\path node at (12.3,6.3) {$3$}; \path node at (13.3,6.3) {$0$}; \path node at (14.3,6.3) {$4$}; \path node at (15.3,6.3) {$1$}; \path node at (16.3,6.3) {$5$};
\path node at (17.3,6.3) {$2$}; \path node at (18.3,6.3) {$6$}; \path node at (19.3,6.3) {$3$}; \path node at (20.3,6.3) {$0$};

\path node at (12.3,7.3) {$5$}; \path node at (13.3,7.3) {$2$}; \path node at (14.3,7.3) {$6$}; \path node at (15.3,7.3) {$3$}; \path node at (16.33,7.3) {$0$};
\path node at (17.3,7.3) {$4$}; \path node at (18.3,7.3) {$1$}; \path node at (19.3,7.3) {$5$}; \path node at (20.3,7.3) {$2$};

\path node at (12,0) {$0$}; \path node at (13,0) {$1$}; \path node at (14,0) {$2$}; \path node at (15,0) {$3$}; \path node at (16,0) {$4$};
\path node at (17,0) {$5$}; \path node at (18,0) {$6$}; \path node at (19,0) {$7$};  \path node at (20,0) {$8$};

\path node at (11,1) {$0$}; \path node at (11,2) {$1$}; \path node at (11,3) {$2$}; \path node at (11,4) {$3$}; \path node at (11,5) {$4$};
\path node at (11,6) {$5$}; \path node at (11,7) {$6$};

\path node at (16,-1) {$d)$};

\filldraw (1,10) circle (2pt); \filldraw (2,10) circle (2pt); \filldraw (3,10) circle (2pt); \filldraw (4,10) circle (2pt);
\filldraw (5,10) circle (2pt); \filldraw (6,10) circle (2pt); \filldraw (7,10) circle (2pt); \filldraw (8,10) circle (2pt); \filldraw (9,10) circle (2pt);
\filldraw (1,11) circle (2pt); \filldraw (2,11) circle (2pt); \filldraw (3,11) circle (2pt); \filldraw (4,11) circle (2pt);
\filldraw (5,11) circle (2pt); \filldraw (6,11) circle (2pt); \filldraw (7,11) circle (2pt); \filldraw (8,11) circle (2pt); \filldraw (9,11) circle (2pt);
\filldraw (1,12) circle (2pt); \filldraw (2,12) circle (2pt); \filldraw (3,12) circle (2pt); \filldraw (4,12) circle (2pt); 
\filldraw (5,12) circle (2pt); \filldraw (6,12) circle (2pt); \filldraw (7,12) circle (2pt); \filldraw (8,12) circle (2pt); \filldraw (9,12) circle (2pt);
\filldraw (1,13) circle (2pt); \filldraw (2,13) circle (2pt); \filldraw (3,13) circle (2pt); \filldraw (4,13) circle (2pt);
\filldraw (5,13) circle (2pt); \filldraw (6,13) circle (2pt); \filldraw (7,13) circle (2pt); \filldraw (8,13) circle (2pt); \filldraw (9,13) circle (2pt);
\filldraw (1,14) circle (2pt); \filldraw (2,14) circle (2pt); \filldraw (3,14) circle (2pt); \filldraw (4,14) circle (2pt);
\filldraw (5,14) circle (2pt); \filldraw (6,14) circle (2pt); \filldraw (7,14) circle (2pt); \filldraw (8,14) circle (2pt); \filldraw (9,14) circle (2pt); 
\filldraw (1,15) circle (2pt); \filldraw (2,15) circle (2pt); \filldraw (3,15) circle (2pt); \filldraw (4,15) circle (2pt);
\filldraw (5,15) circle (2pt); \filldraw (6,15) circle (2pt); \filldraw (7,15) circle (2pt); \filldraw (8,15) circle (2pt); \filldraw (9,15) circle (2pt);
\filldraw (1,16) circle (2pt); \filldraw (2,16) circle (2pt); \filldraw (3,16) circle (2pt); \filldraw (4,16) circle (2pt);
\filldraw (5,16) circle (2pt); \filldraw (6,16) circle (2pt); \filldraw (7,16) circle (2pt); \filldraw (8,16) circle (2pt); \filldraw (9,16) circle (2pt);

 \draw  (1,10)--(9,10); \draw  (1,11)--(9,11);\draw  (1,12)--(9,12);\draw  (1,13)--(9,13);\draw  (1,14)--(9,14);\draw  (1,15)--(9,15);\draw  (1,16)--(9,16);

\draw  (1,10)--(1,16); \draw  (2,10)--(2,16);\draw  (3,10)--(3,16);\draw  (4,10)--(4,16);\draw  (5,10)--(5,16);\draw  (6,10)--(6,16);\draw  (7,10)--(7,16);
\draw  (8,10)--(8,16);\draw  (9,10)--(9,16);

\path node at (1.3,10.3) {$0$}; \path node at (2.3,10.3) {$2$}; \path node at (3.3,10.3) {$4$}; \path node at (4.3,10.3) {$6$}; \path node at (5.3,10.3) {$1$};
\path node at (6.3,10.3) {$3$}; \path node at (7.3,10.3) {$5$}; \path node at (8.3,10.3) {$0$}; \path node at (9.3,10.3) {$2$};

\path node at (1.3,11.3) {$3$}; \path node at (2.3,11.3) {$5$}; \path node at (3.3,11.3) {$0$}; \path node at (4.3,11.3) {$2$}; \path node at (5.3,11.3) {$4$};
\path node at (6.3,11.3) {$6$}; \path node at (7.3,11.3) {$1$}; \path node at (8.3,11.3) {$3$}; \path node at (9.3,11.3) {$5$};

\path node at (1.3,12.3) {$6$}; \path node at (2.3,12.3) {$1$}; \path node at (3.3,12.3) {$3$}; \path node at (4.3,12.3) {$5$}; \path node at (5.3,12.3) {$0$};
\path node at (6.3,12.3) {$2$}; \path node at (7.3,12.3) {$4$}; \path node at (8.3,12.3) {$6$}; \path node at (9.3,12.3) {$1$};

\path node at (1.3,13.3) {$2$}; \path node at (2.3,13.3) {$4$}; \path node at (3.3,13.3) {$6$}; \path node at (4.3,13.3) {$1$}; \path node at (5.3,13.3) {$3$};
\path node at (6.3,13.3) {$5$}; \path node at (7.3,13.3) {$0$}; \path node at (8.3,13.3) {$2$}; \path node at (9.3,13.3) {$4$};

\path node at (1.3,14.3) {$5$}; \path node at (2.3,14.3) {$0$}; \path node at (3.3,14.3) {$2$}; \path node at (4.3,14.3) {$4$}; \path node at (5.3,14.3) {$6$};
\path node at (6.3,14.3) {$1$}; \path node at (7.3,14.3) {$3$}; \path node at (8.3,14.3) {$5$}; \path node at (9.3,14.3) {$0$};
 
\path node at (1.3,15.3) {$1$}; \path node at (2.3,15.3) {$3$}; \path node at (3.3,15.3) {$5$}; \path node at (4.3,15.3) {$0$}; \path node at (5.3,15.3) {$2$};
\path node at (6.3,15.3) {$4$}; \path node at (7.3,15.3) {$6$}; \path node at (8.3,15.3) {$1$}; \path node at (9.3,15.3) {$3$};

\path node at (1.3,16.3) {$4$}; \path node at (2.3,16.3) {$6$}; \path node at (3.3,16.3) {$1$}; \path node at (4.3,16.3) {$3$}; \path node at (5.3,16.3) {$5$};
\path node at (6.3,16.3) {$0$}; \path node at (7.3,16.3) {$2$}; \path node at (8.3,16.3) {$4$}; \path node at (9.3,16.3) {$6$};

\path node at (1,9) {$0$}; \path node at (2,9) {$1$}; \path node at (3,9) {$2$}; \path node at (4,9) {$3$}; \path node at (5,9) {$4$};
\path node at (6,9) {$5$}; \path node at (7,9) {$6$}; \path node at (8,9) {$7$};  \path node at (9,9) {$8$};

\path node at (0,10) {$0$}; \path node at (0,11) {$1$}; \path node at (0,12) {$2$}; \path node at (0,13) {$3$}; \path node at (0,14) {$4$};
\path node at (0,15) {$5$}; \path node at (0,16) {$6$};

\path node at (5,8) {$a)$};

\filldraw (12,10) circle (2pt); \filldraw (13,10) circle (2pt); \filldraw (14,10) circle (2pt); \filldraw (15,10) circle (2pt);
\filldraw (16,10) circle (2pt); \filldraw (17,10) circle (2pt); \filldraw (18,10) circle (2pt); \filldraw (19,10) circle (2pt); \filldraw (20,10) circle (2pt);
\filldraw (12,11) circle (2pt); \filldraw (13,11) circle (2pt); \filldraw (14,11) circle (2pt); \filldraw (15,11) circle (2pt);
\filldraw (16,11) circle (2pt); \filldraw (17,11) circle (2pt); \filldraw (18,11) circle (2pt); \filldraw (19,11) circle (2pt); \filldraw (20,11) circle (2pt);
\filldraw (12,12) circle (2pt); \filldraw (13,12) circle (2pt); \filldraw (14,12) circle (2pt); \filldraw (15,12) circle (2pt); 
\filldraw (16,12) circle (2pt); \filldraw (17,12) circle (2pt); \filldraw (18,12) circle (2pt); \filldraw (19,12) circle (2pt); \filldraw (20,12) circle (2pt);
\filldraw (12,13) circle (2pt); \filldraw (13,13) circle (2pt); \filldraw (14,13) circle (2pt); \filldraw (15,13) circle (2pt);
\filldraw (16,13) circle (2pt); \filldraw (17,13) circle (2pt); \filldraw (18,13) circle (2pt); \filldraw (19,13) circle (2pt); \filldraw (20,13) circle (2pt);
\filldraw (12,14) circle (2pt); \filldraw (13,14) circle (2pt); \filldraw (14,14) circle (2pt); \filldraw (15,14) circle (2pt);
\filldraw (16,14) circle (2pt); \filldraw (17,14) circle (2pt); \filldraw (18,14) circle (2pt); \filldraw (19,14) circle (2pt); \filldraw (20,14) circle (2pt); 
\filldraw (12,15) circle (2pt); \filldraw (13,15) circle (2pt); \filldraw (14,15) circle (2pt); \filldraw (15,15) circle (2pt);
\filldraw (16,15) circle (2pt); \filldraw (17,15) circle (2pt); \filldraw (18,15) circle (2pt); \filldraw (19,15) circle (2pt); \filldraw (20,15) circle (2pt);
\filldraw (12,16) circle (2pt); \filldraw (13,16) circle (2pt); \filldraw (14,16) circle (2pt); \filldraw (15,16) circle (2pt);
\filldraw (16,16) circle (2pt); \filldraw (17,16) circle (2pt); \filldraw (18,16) circle (2pt); \filldraw (19,16) circle (2pt); \filldraw (20,16) circle (2pt);

\draw  (12,10)--(20,10); \draw  (12,11)--(20,11);\draw  (12,12)--(20,12);\draw  (12,13)--(20,13);\draw  (12,14)--(20,14);\draw  (12,15)--(20,15);\draw  (12,16)--(20,16);

\draw  (12,10)--(12,16); \draw  (13,10)--(13,16);\draw  (14,10)--(14,16);\draw  (15,10)--(15,16);\draw  (16,10)--(16,16);\draw  (17,10)--(17,16);\draw  (18,10)--(18,16);\draw  (19,10)--(19,16);\draw  (20,10)--(20,16);

\path node at (12.3,10.3) {$0$}; \path node at (13.3,10.3) {$2$}; \path node at (14.3,10.3) {$4$}; \path node at (15.3,10.3) {$6$}; \path node at (16.3,10.3) {$1$};
\path node at (17.3,10.3) {$3$}; \path node at (18.3,10.3) {$5$}; \path node at (19.3,10.3) {$0$}; \path node at (20.3,10.3) {$2$};

\path node at (12.3,11.3) {$4$}; \path node at (13.3,11.3) {$6$}; \path node at (14.3,11.3) {$1$}; \path node at (15.3,11.3) {$3$}; \path node at (16.3,11.3) {$5$};
\path node at (17.3,11.3) {$0$}; \path node at (18.3,11.3) {$2$}; \path node at (19.3,11.3) {$4$}; \path node at (20.3,11.3) {$6$};

\path node at (12.3,12.3) {$1$}; \path node at (13.3,12.3) {$3$}; \path node at (14.3,12.3) {$5$}; \path node at (15.3,12.3) {$0$}; \path node at (16.3,12.3) {$2$};
\path node at (17.3,12.3) {$4$}; \path node at (18.3,12.3) {$6$}; \path node at (19.3,12.3) {$1$}; \path node at (20.3,12.3) {$3$};

\path node at (12.3,13.3) {$5$}; \path node at (13.3,13.3) {$0$}; \path node at (14.3,13.3) {$2$}; \path node at (15.3,13.3) {$4$}; \path node at (16.3,13.3) {$6$};
\path node at (17.3,13.3) {$1$}; \path node at (18.3,13.3) {$3$}; \path node at (19.3,13.3) {$5$}; \path node at (20.3,13.3) {$0$};

\path node at (12.3,14.3) {$2$}; \path node at (13.3,14.3) {$4$}; \path node at (14.3,14.3) {$6$}; \path node at (15.3,14.3) {$1$}; \path node at (16.3,14.3) {$3$};
\path node at (17.3,14.3) {$5$}; \path node at (18.3,14.3) {$0$}; \path node at (19.3,14.3) {$2$}; \path node at (20.3,14.3) {$4$};
 
\path node at (12.3,15.3) {$6$}; \path node at (13.3,15.3) {$1$}; \path node at (14.3,15.3) {$3$}; \path node at (15.3,15.3) {$5$}; \path node at (16.3,15.3) {$0$};
\path node at (17.3,15.3) {$2$}; \path node at (18.3,15.3) {$4$}; \path node at (19.3,15.3) {$6$}; \path node at (20.3,15.3) {$1$};

\path node at (12.3,16.3) {$3$}; \path node at (13.3,16.3) {$5$}; \path node at (14.3,16.3) {$0$}; \path node at (15.3,16.3) {$2$}; \path node at (16.3,16.3) {$4$};
\path node at (17.3,16.3) {$6$}; \path node at (18.3,16.3) {$1$}; \path node at (19.3,16.3) {$3$}; \path node at (20.3,16.3) {$5$};

\path node at (12,9) {$0$}; \path node at (13,9) {$1$}; \path node at (14,9) {$2$}; \path node at (15,9) {$3$}; \path node at (16,9) {$4$};
\path node at (17,9) {$5$}; \path node at (18,9) {$6$}; \path node at (19,9) {$7$};  \path node at (20,9) {$8$};

\path node at (11,10) {$0$}; \path node at (11,11) {$1$}; \path node at (11,12) {$2$}; \path node at (11,13) {$3$}; \path node at (11,14) {$4$};
\path node at (11,15) {$5$}; \path node at (11,16) {$6$};

\path node at (16,8) {$b)$};	}	
\end{tikzpicture}
\caption{Four  $L(2,1)$-labelings of $P_9\Box P_7$ that determined $L(2,1)$-labelings of Cartesian graph bundle $C_9\Box^{\sigma_\ell}C_7$ according to the cyclic $\ell$-shift $\sigma_\ell$:  \\
a) $\ell =6, f_1(i,j)=\left[2i+3j\right]\mod 7$,\\
b) $\ell =1, f_2(i,j)=\left[2i+4j)\right]\mod 7$, \\
c) $\ell=3, g_1(i,j)=\left[3i+2j\right]\mod 7$, \\
d) $\ell=4, g_2(i,j)=\left[4i+2j\right]\mod 7$.}
\label{primer2}
\end{center}
\end{figure}

{\bf Funding:} This research received no external funding.

{\bf Data Availability Statement:} All relevant data are within the manuscript.

{\bf Conflicts of Interest:} The author declares no conflicts of interest.

{\bf Acknowledgments:} The authors wish to sincerely thank three anonymous reviewers for careful
reading of the first version of the paper and for providing constructive remarks that helped us to
considerably improve the paper.



\end{document}